\documentclass{amsart}

\newtheorem{thm}{Theorem}[section]
\newtheorem{cor}[thm]{Corollary}

\newtheorem{prop}[thm]{Proposition}
\theoremstyle{definition}
\newtheorem{defn}[thm]{Definition}
\theoremstyle{remark}
\newtheorem{rem}[thm]{Remark}
\newtheorem{ex}[thm]{Example}

\begin{document}
\title[Paraquaternionic CR-submanifolds]{\bf Paraquaternionic CR-submanifolds of paraquaternionic K\"{a}hler manifolds and semi-Riemannian submersions}

\author[S. Ianu\c{s}, S. Marchiafava, G.E. V\^{\i}lcu]{Stere Ianu\c{s}, Stefano Marchiafava, Gabriel Eduard V\^{\i}lcu}

\date{}
\maketitle

\abstract In this paper we introduce paraquaternionic
CR-submanifolds of almost paraquaternionic hermitian manifolds and
state some basic results on their differential geometry. We also
study a class of semi-Riemannian submersions from paraquaternionic
CR-submanifolds of paraquaternionic K\"{a}hler manifolds.\\
{\em AMS Mathematics Subject Classification:} 53C15.\\
{\em Key Words and Phrases:} paraquaternionic K\"{a}hler manifold,
foliation, CR-submanifold, semi-Riemannian submersion.

\endabstract

\section{Introduction}

The notion of CR-submanifold of a K\"{a}hler manifold was introduced
by Bejancu \cite{BJC} as a generalization both of totally real and
of holomorphic submanifolds of K\"{a}hler manifolds. Then many
papers appeared studying geometry of CR-submanifolds of K\"{a}hler
manifolds and this notion was further extended to other ambient
spaces; in the monographs \cite{BEJN,CHN3,YK2} we can find the most
significant results concerning CR-submanifolds. A class of examples
of CR-submanifolds of manifolds endowed with different geometric
structures is given in \cite{ORN3}.

On the other hand, the paraquaternionic structures, firstly named
quaternionic structures of second kind, have been introduced by P.
Libermann in \cite{LBM}. The differential geometry of manifolds
endowed with this kind of structures is a very interesting subject
and these manifolds have been intensively studied by many authors
(see, e.g., \cite{AC,AK,BL,DJS,GM,IZ,IMZ,ZMK}). The study of
submanifolds of a paraquaternionic K\"{a}hler manifold is also of
interest and several types of such submanifolds we can find in the
recent literature: paraquaternionic submanifolds \cite{VAC},
K\"{a}hler and para-K\"{a}hler submanifolds \cite{AC,MAR}, normal
semi-invariant submanifolds \cite{AB,BEJ}, lightlike submanifolds
\cite{IMV,IMV2}, $F$-invariant submanifolds \cite{VLC}. In this note
we define a new class of submanifolds of paraquaternionic K\"{a}hler
manifolds, which we call paraquaternionic CR-submanifolds, as a
natural extension of CR-submanifolds in paraquaternionic setting.

The paper is organized as follows: in Section 2 we collect basic
definitions, some formulas and results for later use. In Section 3
we introduce the concept of paraquaternionic CR-submanifold and show
that on the normal bundle of a paraquaternionic CR-submanifold of a
paraquaternionic K\"{a}hler manifold exist two complementary
orthogonal distributions. In Section 4 we investigate the
integrability of the distributions involved in the definition of a
paraquaternionic CR-submanifold. In Section 5, following the same
techniques as in \cite{IIV}, we study the canonical foliation
induced on a paraquaternionic CR-submanifold; conditions on total
geodesicity are derived. We also obtain necessary and sufficient
conditions for a paraquaternionic CR-submanifold of a
paraquaternionic K\"{a}hler manifold to be a ruled submanifold with
respect to the canonical foliation. In Section 6 we define the
paraquaternionic CR-submersions (in the sense of Kobayashi
\cite{KOB}) as semi-Riemannian submersions from paraquaternionic
CR-submanifolds onto an almost paraquaternionic hermitian manifold
and obtain some properties concerning their geometry. In Section 7
we discuss curvature properties of fibers and base manifold for
paraquaternionic CR-submersions.

\section{Preliminaries}

Let $M$ be a smooth manifold endowed with \emph{an almost
paraquaternionic structure} $\sigma$, that is, a rank-3 subbundle
$\sigma$ of $End(TM)$ which admits a local basis
$\lbrace{J_1,J_2,J_3}\rbrace$ on a coordinate neighborhood $U\subset
M$ such that we have:
\begin{equation}\label{1}
J_\alpha^2=-\epsilon_\alpha Id,\
J_{\alpha}J_{\alpha+1}=-J_{\alpha+1}J_{\alpha}=\epsilon_{\alpha+2}J_{\alpha+2}
\end{equation}
where $\epsilon_1=1,\ \epsilon_2=\epsilon_3=-1$ and the indices are
taken from $\{1,2,3\}$ modulo 3. Moreover, the pair $(M,\sigma)$ is
said to be \emph{an almost paraquaternionic manifold} and
$\lbrace{J_1,J_2,J_3}\rbrace$ is called \emph{a canonical local
basis} of $M$.

In an almost paraquaternionic manifold $(M,\sigma)$ we take
intersecting coordinate neighborhoods $U$ and $U'$. Let
$\lbrace{J_1,J_2,J_3}\rbrace$ and $\lbrace{J'_1,J'_2,J'_3}\rbrace$
be canonical local bases of $\sigma$ in $U$ and $U'$ respectively.
Then $\lbrace{J'_1,J'_2,J'_3}\rbrace$ are linear combinations of
$\lbrace{J_1,J_2,J_3}\rbrace$ in $U\cap U'$:
\begin{equation}\label{2bis}
         J'_\alpha=\sum_{\beta=1}^{3}a_{\alpha\beta}J_\beta,\
         \alpha=1,2,3,
         \end{equation}
where $a_{\alpha \beta}$ are functions in $U\cap U'$, %
$\alpha,\beta=1,2,3$ and $A=(a_{{\alpha
\beta}})_{{\alpha,\beta}=1,2,3}\in
       SO(2,1)$.

Let $(M,g)$ be a semi-Riemannian manifold and let $\sigma$ be an
almost paraquaternionic structure on $M$. The metric $g$ is said to
be \emph{adapted} to the almost paraquaternionic structure $\sigma$
if it satisfies:
\begin{equation}\label{2}
g(J_\alpha X,J_\alpha Y)=\epsilon_\alpha g(X,Y),\
\alpha\in\{1,2,3\},
\end{equation}
for all vector fields $X,Y$ on $M$ and any local basis
$\lbrace{J_1,J_2,J_3}\rbrace$ of $\sigma$; or, equivalently, if all
endomorphisms of $\sigma$ are skew-symmetric with respect to $g$. In
this case, $(M,\sigma,g)$ is said to be \emph{an almost
paraquaternionic hermitian manifold}. The existence of
paraquaternionic hermitian structures on manifolds and tangent
bundles has been recently investigated in \cite{IV}.

It is easy to see that any almost paraquaternionic hermitian
manifold is of dimension $4m,\ m\geq 1$, and any adapted metric is
necessarily of neutral signature $(2m,2m)$.

Let $\lbrace{J_1,J_2,J_3}\rbrace$  be a canonical local base of
$\sigma$ in a coordinate neighborhood $U$ of an almost
paraquaternionic hermitian manifold $(M,\sigma,g)$. If we denote by:
         \[
         \Omega_{J_\alpha}(X,Y)=g(X,J_{\alpha}Y),\alpha=1,2,3
         \]
for any vector fields $X$ and $Y$, then, by means of (\ref{2bis}),
we see that:
         \[
         \Omega=\Omega_{J_1}\wedge\Omega_{J_1}-\Omega_{J_2}\wedge\Omega_{J_2}-\Omega_{J_3}\wedge\Omega_{J_3}
         \]
is a globally well-defined 4-form on $M$, called \emph{the
fundamental 4-form} of the manifold.

If $(M,\sigma,g)$ is an almost paraquaternionic hermitian manifold
such that the bundle $\sigma$ is parallel with respect to the
Levi-Civita connection $\overline{\nabla}$ of $g$, then
$(M,\sigma,g)$ is said to be a paraquaternionic K\"{a}hler manifold.
Equivalently, locally defined 1-forms $\omega_1,\omega_2,\omega_3$
exist such that we have for all $ \alpha\in\{1,2,3\}$:
    \begin{equation}\label{3}
    \overline{\nabla}_XJ_{\alpha}=-\epsilon_{\alpha}[\omega_{\alpha+2}(X)J_{\alpha+1}-
    \omega_{\alpha+1}(X)J_{\alpha+2}]
    \end{equation}
for any vector field $X$ on $M$, where the indices are taken from
$\{1,2,3\}$ modulo 3 (see \cite{GM}). Moreover, it can be proved
that an almost paraquaternionic hermitian manifold of dimension
strictly greater than $4$ is paraquaternionic K\"{a}hler if and only
if $\nabla\Omega=0$ (see \cite{IONV}).

We remark that any paraquaternionic K\"{a}hler manifold is an
Einstein manifold, provided that $dim M>4$ (see \cite{BL,GM,IZ}).

Let $(M,g)$ be a semi-Riemannian manifold and let $N$ be an immersed
submanifold of $M$. Then $N$ is said to be a non-degenerate
submanifold of $M$ if the restriction of the semi-Riemannian metric
$g$ to $TN$ is non-degenerate at each point of $N$. We denote by the
same symbol $g$ the semi-Riemannian metric induced by $g$ on $N$ and
by $TN^\perp$ the normal bundle to $N$.

For the rest of this section we will assume that the induced metric
on $N$ is non-degenerate.

Then we have the following orthogonal decomposition:
$$TM=TN\oplus TN^\perp.$$
Also, we denote by $\overline{\nabla}$ and $\nabla$ the Levi-Civita
connection on $M$ and $N$, respectively. Then the Gauss formula is
given by:
\begin{equation}\label{4}
\overline{\nabla}_X Y=\nabla_X Y+B(X,Y)
\end{equation}
for any $X,Y\in\Gamma(TN)$, where
$B:\Gamma(TN)\times\Gamma(TN)\rightarrow\Gamma(TN^\perp)$ is the
second fundamental form of $N$ in $M$.

On the other hand, the Weingarten formula is given by:
\begin{equation}\label{5}
\overline{\nabla}_X \xi=-A_\xi X+\nabla^\perp_X\xi
\end{equation}
for any $X\in\Gamma(TN)$ and $\xi\in\Gamma(TN^\perp)$, where $-A_\xi
X$ is the tangent part of $\overline{\nabla}_X \xi$ and
$\nabla^\perp_X\xi$ is the normal part of $\overline{\nabla}_X \xi$;
$A_\xi$ and $\nabla^\perp$ are called the shape operator of $N$ with
respect to $\xi$ and the normal connection, respectively. Moreover,
$B$ and $A_\xi$ are related by:
\begin{equation}\label{6}
g(B(X,Y),\xi)=g(A_\xi X,Y)
\end{equation}
for any $X,Y\in\Gamma(TN)$ and $\xi\in\Gamma(TN^\perp)$ (see
\cite{ON2}).

If we denote by $\overline{R}$ and $R$ the curvature tensor fields
of $\overline{\nabla}$ and $\nabla$ we have the Gauss equation:
\begin{eqnarray}\label{6b}
&&\overline{R}(X,Y,Z,U)=R(X,Y,Z,U)-g(B(X,Z),B(Y,U))
       +g(B(Y,Z),B(X,U)),
       \end{eqnarray}
for all $X,Y,Z,U\in \Gamma(TN)$.

\section{Paraquaternionic CR-submanifolds}

\begin{defn}\label{3.1}
Let $N$ be an $n$-dimensional non-degenerate submanifold of an
almost paraquaternionic hermitian manifold $(M,\sigma,g)$. We say
that $(N,g)$ is \emph{a paraquaternionic CR-submanifold} of $M$ if
there exists a non-degenerate distribution
$\mathcal{D}:x\rightarrow\mathcal{D}_x\subseteq T_xN$ such that on any $U\cap N$ we have:\\
i. $\mathcal{D}$ is a paraquaternionic distribution, i.e.
\begin{equation}\label{7}
J_\alpha\mathcal{D}_x=\mathcal{D}_x,\ \alpha\in\{1,2,3\}
\end{equation}
and\\
ii. $\mathcal{D}^\perp$ is a totally real distribution, i.e.
\begin{equation}\label{8}
J_\alpha\mathcal{D}^\perp_x\subset T_x^\perp N,\ \alpha\in\{1,2,3\}
\end{equation}
for any local basis $\lbrace{J_1,J_2,J_3}\rbrace$ of $\sigma$ on $U$
and $x\in U\cap N$, where $\mathcal{D}^\perp$ is the orthogonal
complementary distribution to $\mathcal{D}$ in $TN$.
\end{defn}
\begin{defn}
A non-degenerate submanifold $N$ of an almost paraquaternionic
hermitian manifold  $(M,\sigma,g)$ is called \emph{a
paraquaternionic} (respectively, \emph{totally real}) submanifold if
$\mathcal{D}^\perp={0}$ (respectively, $\mathcal{D}={0}$). A
paraquaternionic CR-submanifold is said to be \emph{proper} if it is
neither paraquaternionic nor totally real.
\end{defn}
\begin{ex}
i. The canonical immersion of $P^n\mathbb{B}(c)$ into
$P^m\mathbb{B}(c)$, where $n\leq m$, provides us a very natural
example of paraquaternionic submanifold (see \cite{MAR}).\\
ii. The real projective space $P^n_s\mathbb{R}(\frac{c}{4})$ is a
totally-real submanifold of the paraquaternionic projective space
$P^n\mathbb{B}(c)$, where $s\in\{0,...,n\}$ denotes the index of the
manifold, defined as the dimension of the largest negative definite
vector subspace of the tangent space. \\
iii. Let $(M_1,g_1,\sigma_1)$ and $(M_2,g_2,\sigma_2)$ be two
paraquaternionic K\"{a}hler manifolds. If $U_1$ and $U_2$ are open
subsets of $M_1$ and $M_2$ respectively, on which local basis
$\{J^{(1)}_{1},J^{(1)}_{2},J^{(1)}_{3}\}$  and
$\{J^{(2)}_{1},J^{(2)}_{2},J^{(2)}_{3}\}$ for $\sigma_{1}$ and
$\sigma_{2}$ respectively, are defined, then the product manifold
$U=U_1\times U_2$ can be endowed with an almost paraquaternionic
hermitian non-K\"{a}hler structure $(g,\sigma)$ (see \cite{VLC}).
Now, if $N_1$ is a paraquaternionic submanifold of $U_1$ and $N_2$
is a totally-real submanifold of $U_2$, then $N=N_1\times N_2$ is a
proper paraquaternionic $CR$-submanifold of the almost
paraquaternionic hermitian
manifold $(U,g,\sigma)$.\\
iv. A large class of examples of proper paraquaternionic
CR-submanifolds can be constructed using the paraquaternionic
momentum map \cite{VKM} and the technique from \cite{ORN3}. Suppose
that a Lie group $G$ acts freely and isometrically on the
paraquaternionic K\"{a}hler manifold $(M,\sigma,g)$, preserving the
fundamental 4-form $\Omega$ of the manifold. We denote by
$\mathfrak{g}$ the Lie algebra of $G$, by $\mathfrak{g}^*$ its dual
and by $V$ the unique Killing vector field corresponding to a vector
$V^*\in\mathfrak{g}$. Then there exists a unique section $f$ of
bundle $\mathfrak{g}^*\otimes\sigma$ such that (see \cite{VKM})
\begin{equation}
\nabla f_{V^*}=\theta_{V^*},
\end{equation}
for all $V^*\in\mathfrak{g}$, where the section $\theta_{V^*}$ of
the bundle $\Omega^1(\sigma)$ with values in $\sigma$ is well
defined globally by
\[
\theta_{V^*}(X)=\sum_{\alpha=1}^{3}\omega_\alpha(V,X)J_\alpha,\
\forall X\in TM.
\]
Moreover, the group $G$ acts by isometries on the pre-image
$f^{-1}(0)$ of the zero-section $0\in\mathfrak{g}^*\otimes\sigma$.
Similarly as in \cite{ORN3}, we have the decomposition
\[
T_x(f^{-1}(0))=T_x(G\cdot x)\oplus H_x, \forall x\in M,
\]
where $G\cdot x$ represents the orbit of $G$ through $x$, supposed
to be non-degenerate, and $H_x$ is the orthogonal complementary
subspace of $T_x(G\cdot x)$ in $T_x(f^{-1}(0))$. Because $H_x$ is
invariant under the action of $\sigma$ and $T_x(G\cdot x)$ is
totally real, we can state now the following result.
\begin{prop}\label{example}
If $f^{-1}(0)$ is a smooth submanifold of a paraquaternionic
K\"{a}hler manifold $(M,\sigma,g)$, then $f^{-1}(0)$ is a proper
paraquaternionic $CR$-submanifold of $M$.
\end{prop}

We remark that, in general, $f^{-1}(0)$ is not a differentiable
submanifold of $M$, but always we can take a subset $N\subset
f^{-1}(0)$ which is invariant under the action of $G$ and which is a
submanifold of $M$. A particular example is given in \cite{BV}, as a
paraquaternionic version of the example constructed by Galicki and
Lawson in \cite{GL}: if $p$ and $q$ are distinct and relatively
prime natural numbers, then we have the action of the Lie group
$G=\{e^{jt}|t\in\mathbb{R}\}$ on $P^2\mathbb{B}$ defined by
\[
\phi_{p,q}(t)\cdot[u_0,u_1,u_2]:=[e^{jqt}u_0,e^{jpt}u_1,e^{jpt}u_2],
\]
where $e^{jt}={\rm cosh}t+j{\rm sinh}t$ and $[u_0,u_1,u_2]$ are
homogenous coordinates on $P^2\mathbb{B}$. We can see that this
action is free, isometric and preserves the para-quaternionic
structure on $P^2\mathbb{B}$ and, moreover, we have that the
pre-image by the momentum map $f_{p,q}:P^2\mathbb{B}\rightarrow{\rm
Im}\mathbb{B}$ of the zero-section $0\in{\rm Im}\mathbb{B}$ is (see
also \cite{VKM}):
\[
f_{p,q}^{-1}(0)=\{[u_0,u_1,u_2]\in
P^2\mathbb{B}|q\bar{u}_0ju_0+p\bar{u}_1ju_1+p\bar{u}_2ju_2=0\}.
\]
Finally, we conclude that the subset $N$ of the regular points of
$f_{p,q}^{-1}(0)$, given by
\[
N=\{[u_0,u_1,u_2]\in
f_{p,q}^{-1}(0)|q^2|u_0|^2+p^2|u_1|^2+p^2|u_2|^2\neq 0\}
\]
is a a proper paraquaternionic $CR$-submanifold of $P^2\mathbb{B}$.
\end{ex}
\begin{defn}
Let $N$ be a paraquaternionic CR-submanifold of an
almost paraquaternionic hermitian manifold  $(M,\sigma,g)$. Then we say that:\\
i. $N$ is $\mathcal{D}$-geodesic if $B(X,Y)=0,\ \forall X,Y\in\Gamma(\mathcal{D})$;\\
ii. $N$ is $\mathcal{D}^\perp$-geodesic if $B(X,Y)=0,\ \forall
X,Y\in\Gamma(\mathcal{D}^\perp)$;\\
iii. $N$ is mixed geodesic if $B(X,Y)=0,\ \forall
X\in\Gamma(\mathcal{D}), Y\in\Gamma(\mathcal{D}^\perp)$;\\
iv. $N$ is mixed foliated if $N$ is mixed geodesic and $\mathcal{D}$
is integrable.
\end{defn}

We may easily prove the next result (see \cite{AM2,VAC}):
\begin{prop}
Any paraquaternionic submanifold of a paraquaternionic K\"{a}hler
manifold is a totally geodesic paraquaternionic K\"{a}hler
submanifold.
\end{prop}

By using this proposition, we deduce the next consequences.

\begin{cor}
Let $(N,g)$ be a paraquaternionic submanifold of a paraquaternionic
K\"{a}hler manifold $(M,\sigma,g)$. Then:\\
i. $dim N=4n, n\geq 1$ and the signature of $g|_{TN}$ is (2n,2n);\\
ii. $N$ is an Einstein manifold, provided that $dim N>4$.
\end{cor}
\begin{cor}
The paraquaternionic submanifolds of $\mathbb{R}^{4m}_{2m}$ and of
paraquaternionic projective space $P^m\mathbb{B}$ are locally
isometric with $\mathbb{R}^{4n}_{2n}$ and $P^n\mathbb{B}$,
respectively, where $n\leq m$.
\end{cor}

Next, let $(N,g)$ be a paraquaternionic CR-submanifold of a
paraquaternionic K\"{a}hler manifold $(M,\sigma,g)$. We put
$\nu_{\alpha x}=J_\alpha(D_x^\perp),\ \alpha\in\{1,2,3\}$, and
$\nu_x^\perp=\nu_{1x}\oplus\nu_{2x}\oplus\nu_{3x}$, and remark that
$\nu_{1x},\ \nu_{2x},\ \nu_{3x}$ are mutually orthogonal
non-degenerate vector subspaces of $T_xN^\perp$, for any $x\in U\cap
N$. We also note that the subspaces $\nu_{\alpha x}$ depends on the
choice of the local base $(J_\alpha)_\alpha$, while $\nu^\bot_x$
does'nt depend from it.

\begin{prop}
Let $(N,g)$ be a paraquaternionic CR-submanifold of a
paraquaternionic K\"{a}hler manifold $(M,\sigma,g)$. Then we have:\\
i. $J_\alpha(\nu_{\alpha x})=\mathcal{D}_x^\perp,\ \forall
x\in U\cap N,\ \alpha\in\{1,2,3\}$;\\
ii. $J_\alpha(\nu_{\beta x})=\nu_{\gamma x}$, for any even
permutation $(\alpha,\beta,\gamma)$ of $(1,2,3)$ and $x\in
U\cap N$;\\
iii. The mapping $\nu^\perp:x\in N\rightarrow\nu_x^\perp$ defines a
non-degenerate distribution of dimension $3s$, where
$s=dimD_x^\perp$;\\
iv. $J_\alpha(\nu_{x})=\nu_x,\ \forall x\in U\cap N,\
\alpha\in\{1,2,3\}$, where $\nu$ is the complementary orthogonal
distribution to $\nu^\perp$ in $TN^\perp$.
\end{prop}
\begin{proof}
The assertions i. and ii. are obvious from (\ref{1}). The assertion
iii. follows from the Definition \ref{3.1}, since $J_\alpha,\
\alpha\in\{1,2,3\}$, are automorphisms of $T_xN$, for any $x\in
U\cap N$. Concerning the proof of (iv.), firstly one can remark that
$\nu_x$ does'nt depend from the choice of the local base
$(J_\alpha)_\alpha$ and secondly, one could observe that the
subspace
$TN_x^Q=\mathcal{D}_x\oplus\mathcal{D}_x^\perp\oplus\nu_x^\perp$ is
the paraquaternionic subspace generated by $TN_x$, hence its
orthogonal $\nu$ is also a paraquaternionic subspace, i.e. closed
under $J_\alpha$, $\alpha\in\{1,2,3\}$. This completes the proof of
the proposition.
\end{proof}

\section{Integrability of Distributions}

\begin{thm}\label{4.1AAA}
The distribution $\mathcal{D}^\perp$ is integrable.
\end{thm}
\begin{proof}
For any $U,V\in\Gamma(\mathcal{D}^\perp)$ and
$X\in\Gamma(\mathcal{D})$ we have:
\begin{eqnarray}
g(\nabla_UV,X)&=&g(\overline{\nabla}_UV,X)\nonumber\\
&=&\epsilon_\alpha g(J_\alpha\overline{\nabla}_UV,J_\alpha X)\nonumber\\
&=&\epsilon_\alpha
g(\epsilon_{\alpha}[\omega_{\alpha+2}(U)J_{\alpha+1}V-
    \omega_{\alpha+1}(U)J_{\alpha+2}V]+\overline{\nabla}_UJ_\alpha V,J_\alpha X)\nonumber\\
&=&\epsilon_\alpha g(\overline{\nabla}_UJ_\alpha V,J_\alpha
X)\nonumber
\end{eqnarray}
and by using (\ref{5}) we obtain:
\begin{equation}\label{11}
g(\nabla_UV,X)=-\epsilon_\alpha g(A_{J_\alpha V}U,J_\alpha X).
\end{equation}

On the other hand, if we take $U,V\in\Gamma(\mathcal{D}^\perp)$ and
$C\in\Gamma(TN)$ we have:
\begin{eqnarray}
g(A_{J_\alpha V}U,C)&=&-g(\overline{\nabla}_CJ_\alpha V,U)\nonumber\\
&=&g(\epsilon_{\alpha}[\omega_{\alpha+2}(C)J_{\alpha+1}V-
    \omega_{\alpha+1}(C)J_{\alpha+2}V],U)-g(J_\alpha \overline{\nabla}_CV,U)\nonumber
\end{eqnarray}
hence
$$g(A_{J_\alpha V}U,C)=-g(J_\alpha\overline{\nabla}_CV,U).$$
Moreover, from $C\cdot g(J_\alpha U,V)=0$, one has
$$g(J_\alpha \overline{\nabla}_C U,V)=g(J_\alpha\overline{\nabla}_CV,U).$$

We conclude:
\begin{equation}\label{12}
g(A_{J_\alpha V}U,C)=g(A_{J_\alpha U}V,C).
\end{equation}

From (\ref{11}) and (\ref{12}) we deduce that for any
$U,V\in\Gamma(\mathcal{D}^\perp)$ and $X\in\Gamma(\mathcal{D})$ we
have:
$$g(\nabla_UV-\nabla_VU,X)=0,$$
which implies $[U,V]\in\Gamma(\mathcal{D}^\perp),\ \forall
U,V\in\Gamma(\mathcal{D}^\perp)$. Thus $\mathcal{D}^\perp$ is
integrable.
\end{proof}

\begin{thm}\label{4.2b}
The paraquaternionic distribution $\mathcal{D}$ is integrable if and
only if $N$ is $\mathcal{D}$-geodesic.
\end{thm}
\begin{proof}
Similarly as in above theorem, we obtain:
$$
g(\nabla_XY,U)=\epsilon_\alpha g(\overline{\nabla}_XJ_\alpha
Y,J_\alpha U)
$$
for any $X,Y\in\Gamma(\mathcal{D})$ and
$U\in\Gamma(\mathcal{D}^\perp)$, and taking into account (\ref{4})
we obtain:
\begin{equation}\label{13}
g(\nabla_XY,U)=\epsilon_\alpha g(B(X,J_\alpha Y),J_\alpha U).
\end{equation}

If we suppose that $N$ is $\mathcal{D}$-geodesic, from (\ref{13}) we
derive $\nabla_XY\in\Gamma(\mathcal{D})$, which implies
$[X,Y]\in\Gamma(\mathcal{D})$. Thus $\mathcal{D}$ is integrable.

Conversely, if we suppose that $\mathcal{D}$ is integrable, then the
leaves are invariant of the paraquaternionic structure and so they
are totally geodesic in $N$. In particular, $N$ is
$\mathcal{D}$-geodesic.
\end{proof}

\begin{rem}
If $\mathcal{Q}$ is a non-degenerate distribution on a
semi-Riemannian manifold $(M,g)$, then we can consider a
well-defined $\mathcal{Q}^\perp$-valued vector field on $N$, called
the mean curvature vector of $\mathcal{Q}$ (see \cite{BJCF}), given
by:
$$H^\mathcal{Q}=\frac{1}{q}\sum_{i=1}^{q}\theta_ih^\mathcal{Q}(E_i,E_i),$$
where $h^\mathcal{Q}$ is the second fundamental forms of
$\mathcal{Q}$, $q=dim\mathcal{Q}$, $\{E_1,...,E_{q}\}$ is a
pseudo-orthonormal basis of $\mathcal{Q}$ and
$\theta_i=g(E_i,E_i)\in\{-1,1\},\ \forall i\in\{1,...,q\}$.

The distribution $\mathcal{Q}$ is said to be \emph{minimal} if the
mean curvature vector $H^\mathcal{Q}$ of $\mathcal{Q}$ vanishes
identically.
\end{rem}

\begin{thm}\label{4.2AAA}
The paraquaternionic distribution $\mathcal{D}$ is minimal.
\end{thm}
\begin{proof}
For any $X\in\Gamma(D)$ and $U\in\Gamma(D^\perp)$ we obtain
similarly as in above theorems:
\begin{equation}\label{14}
       g(\nabla _X X,U)=\epsilon_\alpha g(A_{J_\alpha  U} J_\alpha  X,X)
       \end{equation}
and
\begin{equation}\label{15}
g(\nabla _{J_\alpha  X} J_\alpha  X,U) = - g(A_{J_\alpha U}
X,J_\alpha  X)= - g(A_{J_\alpha U} J_\alpha  X,X),
\end{equation}
for all $\alpha\in\{1,2,3\}$.

From (\ref{14}) and (\ref{15}) we deduce:
$$
 g(\nabla _X X+\epsilon_\alpha\nabla _{J_\alpha  X} J_\alpha  X,U)=0
$$
and so
\begin{equation}\label{16}
 h^\mathcal{D}(X,X)+\epsilon_\alpha h^\mathcal{D}(J_\alpha  X, J_\alpha  X)=0,\
 \forall\alpha\in\{1,2,3\}.
\end{equation}

From (\ref{16}) it follows that
\[
 h^\mathcal{D}(X,X)=- h^\mathcal{D}(J_1X,J_1X)=- h^\mathcal{D}(J_2J_3X,J_2J_3X)=- h^\mathcal{D}(J_3X,J_3X)=- h^\mathcal{D}(X,X)
\]
and so
\begin{equation}\label{17}
h^\mathcal{D}(X,X)=0,\ \forall X\in\Gamma(D).
\end{equation}

We obtain $H^\mathcal{D}=0$ and the assertion follows. We can remark
that the relation (\ref{17}) not imply $h^\mathcal{D}=0$, since
$h^\mathcal{D}$ is not symmetric in general, unless $\mathcal{D}$ is
integrable.
\end{proof}

\begin{rem}
For any paraquaternionic CR-submanifold $(N,g)$ of an almost
paraquaternionic hermitian manifold $(M^{4m},\sigma,g)$, having
$dim\mathcal{D}=4r$ and $dim\mathcal{D}^\perp=p$, we can choose a
local pseudo-orthonormal frame in $M$:
$$\{e_1,,...,e_r,e_{r+1},...,e_{r+p},e_{r+p+1},...,e_m,J_1e_1,...,J_1e_m,J_2e_1,...,J_2e_m,
J_3e_1,...,J_3e_m\}$$ such that restricted to $N$,
$\{e_i,J_1e_i,J_2e_i,J_3e_i\}_{i\in\{1,...,r\}}$ are in
$\mathcal{D}$ and $\{e_{r+1},...,e_{r+p}\}$ are in
$\mathcal{D}^\perp$.

Let $\{\omega^1,...,\omega^{4r}\}$ be the $4r$ 1-forms on $N$
satisfying:
\begin{equation}\label{107}
\omega^i(Z)=0,\ \omega_i(E_j)=\delta_{ij},\ i,j\in\{1,...,4r\},
\end{equation}
for any $Z\in\Gamma(\mathcal{D}^\perp)$, where $E_j=e_j$,
$E_{r+j}=J_1e_j$, $E_{2r+j}=J_2e_j$ and $E_{3r+j}=J_3e_j$,
$j\in\{1,...,r\}$. Then $\omega=\omega^1\wedge...\wedge\omega^{4r}$
does not depend on the particular base $\{J_1, J_2, J_3\}$ and
defines a $4r$-form on $N$. Therefore, we have
\begin{equation}
d\omega=\sum_{i=1}^{4r}(-1)^i\omega^1\wedge...\wedge
d\omega^i...\wedge\omega^{4r}
\end{equation}
and taking account of (\ref{107}), we deduce that $d\omega=0$ if and
only if
\begin{equation}\label{108}
d\omega(Z_1,Z_2,X_1,...,X_{4r-1})=0
\end{equation}
and
\begin{equation}\label{109}
d\omega(Z_1,X_1,...,X_{4r})=0
\end{equation}
for any $Z_1,Z_2\in\Gamma(\mathcal{D}^\perp)$ and
$X_1,...,X_{4r}\in\Gamma(\mathcal{D})$. We remark now easily that
(\ref{108}) holds if and only if $\mathcal{D}^\perp$ is integrable
and (\ref{109}) holds if and only if $\mathcal{D}$ is minimal (see
also \cite{CHN4}). On another hand, if
$\{\omega^{4r+1},...,\omega^{4r+p}\}$ is the dual frame to the
pseudo-orthonormal frame $\{e_{r+1},...,e_{r+p}\}$ of
$\mathcal{D}^\perp$, we can define a $p$-form $\omega^\perp$ on $M$
by $\omega^\perp=\omega^{4r+1}\wedge...\wedge\omega^{4r+p}$.
Similarly we find that $\omega^\perp$ is closed if $\mathcal{D}$ is
integrable and $\mathcal{D}^\perp$ is minimal.

Consequently, from Theorems \ref{4.1AAA} and \ref{4.2AAA}, we obtain
the following result.

\begin{thm}
Let $N$ be a closed paraquaternionic CR-submanifold of a
paraquaternionic K\"{a}hler manifold $(M,\sigma,g)$. Then the
$4r$-form $\omega$ is closed and defines a canonical de Rham
cohomology class $[\omega]$ in $H^{4r}(M,\mathbb{R})$. Moreover,
this cohomology class is non-trivial if $\mathcal{D}$ is integrable
and $\mathcal{D}^\perp$ is minimal.
\end{thm}

\end{rem}

\section{Canonical foliations on paraquaternionic CR-submanifolds}

Since the totally real distribution $\mathcal{D}^\perp$ of a
paraquaternionic CR-submanifold $N$ of a paraquaternion K\"{a}hler
manifold $(M,\sigma,g)$ is always integrable we conclude that we
have a foliation $\mathfrak{F}^\perp$ on $N$ with structural
distribution $\mathcal{D}^\perp$ and transversal distribution
$\mathcal{D}$ (see \cite{BJCF}). We say that $\mathfrak{F}^\perp$ is
the canonical totally real foliation on $N$.

\begin{thm}\label{4.1}
Let $N$ be a paraquaternionic CR-submanifold of a paraquaternionic
K\"{a}hler manifold $(M,\sigma,g)$. The next assertions are
equivalent:
\\
i. $\mathfrak{F}^\perp$ is totally geodesic;
\\
ii. $B(X,Y)\in\Gamma(\nu)$, $\forall X\in\Gamma(\mathcal{D})$,
$Y\in\Gamma(\mathcal{D}^\perp)$;
\\
iii. $A_NX\in\Gamma(\mathcal{D}^\perp)$, $\forall
X\in\Gamma(\mathcal{D}^\perp)$, $N\in\Gamma(\nu^\perp)$;
\\
iv. $A_NY\in\Gamma(\mathcal{D})$, $\forall Y\in\Gamma(\mathcal{D})$,
$N\in\Gamma(\nu^\perp)$.
\end{thm}
\begin{proof}
By using (\ref{2}), (\ref{3}), (\ref{4}) and (\ref{5}) we obtain for
any $X,Z\in\Gamma(\mathcal{D}^\perp)$ and $Y\in\Gamma(\mathcal{D})$:
\begin{eqnarray}
g(J_\alpha(\nabla_XZ),Y)&=&-g(\overline{\nabla}_XZ,J_\alpha Y)\nonumber\\
&=&g(\epsilon_{\alpha}[\omega_{\alpha+2}(X)J_{\alpha+1}Z-
    \omega_{\alpha+1}(X)J_{\alpha+2}Z]+\overline{\nabla}_XJ_\alpha
       Z,Y)\nonumber\\
       &=&g(-A_{J_\alpha Z}X+\nabla_X^{\perp}J_\alpha
Z,Y)\nonumber\\
       &=&-g(A_{J_\alpha Z}X,Y)\nonumber
       \end{eqnarray}
and taking into account (\ref{6}) we derive:
\begin{equation}\label{18}
g(J_\alpha(\nabla_XZ),Y)=-g(B(X,Y),J_\alpha Z).
\end{equation}
i. $\Rightarrow$ ii. If $\mathfrak{F}^\perp$ is totally geodesic,
then $\nabla_XZ\in\Gamma(\mathcal{D}^\perp)$, for
$X,Z\in\Gamma(\mathcal{D}^\perp)$ and from (\ref{18}) we derive:
$$g(B(X,Y),J_\alpha
Z)=0$$ and the implication is clear.
\\
ii. $\Rightarrow$ i. If we suppose $B(X,Y)\in\Gamma(\mu)$, $\forall
X\in\Gamma(\mathcal{D})$, $Y\in\Gamma(\mathcal{D}^\perp)$, then from
(\ref{18}) we derive:
$$g(J_\alpha(\nabla_XZ),Y)=0$$
and we conclude $\nabla_XZ\in\Gamma(\mathcal{D}^\perp)$. Thus
$\mathfrak{F}^\perp$ is totally geodesic.
\\
ii. $\Leftrightarrow$ iii. This equivalence is clear from (\ref{6}).
\\
iii. $\Leftrightarrow$ iv. This equivalence is true because $A_N$ is
a self-adjoint operator.

\end{proof}

\begin{cor}\label{4.2}
If $N$ is a mixed-geodesic paraquaternionic CR-submanifold of a
paraquaternionic K\"{a}hler manifold $(M,\sigma,g)$, then the
canonical totally real foliation $\mathfrak{F}^\perp$ on $N$ is
totally geodesic.
\end{cor}

\begin{cor}\label{4.3}
Let $N$ be a paraquaternionic CR-submanifold of a paraquaternionic
K\"{a}hler manifold $(M,\sigma,g)$ with $\nu=0$. Then the canonical
totally real foliation $\mathfrak{F}^\perp$ on $N$ is totally
geodesic if and only if $N$ is mixed geodesic.
\end{cor}

We note that the condition $\nu=0$ in Corollary \ref{4.3}
characterizes the following interesting class of paraquaternionic
CR-submanifolds (compare with the definition in the complex case
given in \cite{YK2}, pag. 78).

\begin{defn}
A paraquaternionic CR-submanifold $N$ of a paraquaternionic
K\"{a}hler manifold $(M,\sigma,g)$ is called {\it generic} if the
restriction $TM_{|N}$ to $N$ of the tangent bundle of ambient
manifold $M$ is generated over the paraquaternions by the tangent
bundle of $N$, i.e. if $TN+\sigma TN=TM_{|N}$ or, equivalently, if
$\nu=0$.
\end{defn}

\begin{rem}
The proper paraquaternionic CR-submanifolds given by Proposition
\ref{example} are generic because we have
\[T_x^\perp(f^{-1}(0))=\sigma T_x(G\cdot x).\]
\end{rem}

\begin{prop}
Let $N$ be a generic paraquaternionic CR-submanifold of a
paraquaternionic K\"{a}hler manifold $(M,\sigma,g)$. Then the
paraquaternionic distribution $\mathcal{D}$ is geodesic.
\end{prop}
\begin{proof}
$N$ generic means that
\[ TM_{|N} = \mathcal{D}\oplus{\mathcal D}^\bot\oplus \sigma\mathcal{D}^\bot
\]
where
\[
\nu^\bot =\sigma\mathcal{D}^\bot,\ \ TN = \mathcal{D}\oplus{\mathcal
D}^\bot \, .
\]

From (\ref{2}) and (\ref{13}) we have
\[ g(J_\alpha\nabla_XY,J_\alpha U) =  \epsilon_\alpha g(\nabla_XY,U)) = g(B(X,J_\alpha Y),
J_\alpha U),
\]
for all $X,Y \in \Gamma(\mathcal D),\, U \in \Gamma({\mathcal
D}^\bot)$ and since $TN^\bot=\nu^\bot$,
\[
(J_\alpha \nabla_XY)^\bot = B(X,J_\alpha Y)
\]
where $W^\bot$ means the component of the vector $W \in TN$ in the
or\-tho\-go\-nal to $TN$, that is in $\nu^\bot$.

More explicitly, we have
\[
(J_\alpha\nabla^\mathcal{D}_XY
+J_\alpha\nabla^{\mathcal{D}^\bot}_XY)^\bot = B(X,J_\alpha Y)
\]
where $\nabla^\mathcal{D}_XY, \nabla^{\mathcal{D}^\bot}_XY$ are the
$\mathcal{D}, \mathcal{D}^\bot$ components of $\nabla_XY$
respectively.

Equivalently, since $J_\alpha\nabla^\mathcal{D}_XY \in \mathcal{D}$,
we have
\[
J_\alpha\nabla^{\mathcal{D}^\bot}_XY = B(X,J_\alpha Y)
\]
that is
\begin{equation}
\nabla^{\mathcal{D}^\bot}_XY = -\epsilon_\alpha J_\alpha
B(X,J_\alpha Y)
\end{equation}

As a consequence, we have
\[
\epsilon_\alpha J_\alpha B(X,J_\alpha Y) = \epsilon_\beta J_\beta
B(X,J_\beta Y),\ \forall \alpha, \beta
\]
and, for $J_\alpha J_\beta=\epsilon_\gamma J_\gamma$,
\[
B(X,J_\alpha Y) = -\epsilon_\beta \epsilon_\gamma J_\gamma
B(X,J_\beta Y).
\]
Hence
\begin{equation}\label{basic}
B(X,Y) = -\epsilon_\alpha\epsilon_\beta J_\gamma B(X,J_\gamma Y)
\end{equation}
and also
\begin{equation}\label{comparison}
 B(X,J_\gamma Y) = J_\gamma B(X,Y)
\end{equation}
since $\epsilon_\gamma\epsilon_\alpha\epsilon_\beta=1$.

Note that from (\ref{comparison}) and symmetry of $B$ it follows
\[
 B(X,J_\gamma Y) = B(Y,J_\gamma X)
\]
from which we deduce:
\begin{equation}\label{comparison2}
B(X,Y)= -\epsilon_\gamma B(J_\gamma X,J_\gamma Y).
\end{equation}

By applying repeatedly the (\ref{comparison2}) we find
\[
B(X,Y)=-B(J_1X,J_1Y)=-B(J_2J_3X,J_2J_3Y)=-B(J_3X,J_3Y)=-B(X,Y)
\]
hence,
\[
B=0 \,.
\]
\end{proof}

\begin{defn} \cite{BJCF}
A submanifold $N$ of a semi-Riemannian manifold $(M,g)$ is said to
be a \emph{ruled submanifold} if it admits a foliation whose leaves
are totally geodesic submanifolds immersed in $(M,g)$.
\end{defn}

\begin{defn}
A paraquaternionic CR-submanifold of a paraquaternionic K\"{a}hler
manifold which is a ruled submanifold with respect to the foliation
$\mathfrak{F}^\perp$ is called \emph{totally real ruled
paraquaternionic CR-submanifold}.
\end{defn}

\begin{thm}\label{4.5}
Let $N$ be a paraquaternionic CR-submanifold of a paraquaternionic
K\"{a}hler manifold $(M,\sigma,g)$. Then the following assertions
are mutually equivalent:
\\
i. $N$ is a totally real ruled paraquaternionic CR-submanifold.
\\
ii. $N$ is $\mathcal{D}^\perp$-geodesic and:
$$B(X,Y)\in\Gamma(\nu),\ \forall
X\in\Gamma(\mathcal{D}),\ Y\in\Gamma(\mathcal{D}^\perp).$$
iii. The
subbundle $\nu^\perp$ is $\mathcal{D}^\perp$-parallel, i.e:
$$\nabla_X^\perp J_\alpha Z\in\Gamma(\nu^\perp),\ \forall X,Z\in
\Gamma(\mathcal{D}^\perp),\ \alpha\in\{1,2,3\}$$ and the second
fundamental form satisfies:
$$B(X,Y)\in\Gamma(\nu),\ \forall
X\in\Gamma(\mathcal{D}^\perp),\ Y\in\Gamma(TN).$$
iv. The shape
operator satisfies:
$$A_{J_\alpha Z}X=0,\ \forall X,Z\in \Gamma(\mathcal{D}^\perp),\ \alpha\in\{1,2,3\}$$
and
$$A_N X\in\Gamma(\mathcal{D}),\ \forall X\in\Gamma(\mathcal{D}^\perp),\ N\in\Gamma(\nu).$$
\end{thm}
\begin{proof}
i. $\Leftrightarrow$ ii. This equivalence follows from Theorem
\ref{4.1} since for any $X,Z\in\Gamma(\mathcal{D}^\perp)$ we have:
\begin{eqnarray}
       \overline{\nabla}_XZ&=&\nabla_XZ+B(X,Z)\nonumber\\
       &=&\nabla_X^{D^\perp}Z+h^\perp(X,Z)+B(X,Z)\nonumber
       \end{eqnarray}
and thus the leaves of $\mathcal{D}^\perp$ are totally geodesic
immersed in $M$ if and only if $h^\perp=0$ and $N$ is
$\mathcal{D}^\perp$-geodesic.
\\
i. $\Leftrightarrow$ iii. If $X,Z\in\Gamma(\mathcal{D}^\perp)$ and
$N\in\Gamma(\nu)$, then we have:
\begin{eqnarray}
g(\overline{\nabla}_XZ,N)&=&\epsilon_\alpha g(J_\alpha\overline{\nabla}_XZ,J_\alpha N)\nonumber\\
&=&\epsilon_\alpha
g(\epsilon_{\alpha}[\omega_{\alpha+2}(X)J_{\alpha+1}Z-
    \omega_{\alpha+1}(X)J_{\alpha+2}Z]+\overline{\nabla}_XJ_\alpha Z,J_\alpha N)\nonumber\\
       &=&\epsilon_\alpha g(-A_{J_\alpha Z}X+\nabla_X^{\perp}J_\alpha
Z,J_\alpha N)\nonumber
       \end{eqnarray}
and thus we obtain:
\begin{equation}\label{19}
g(\overline{\nabla}_XZ,N)=\epsilon_\alpha g(\nabla^\perp_X J_\alpha
Z,J_\alpha N).
\end{equation}

Similarly we find:
\begin{equation}\label{20}
g(\overline{\nabla}_XZ,U)=-\epsilon_\alpha g(B(X,J_\alpha
U),J_\alpha Z),\ \forall X,Z\in\Gamma(\mathcal{D}^\perp),\
U\in\Gamma(\mathcal{D}).
\end{equation}

On the other hand, from (\ref{4}) we deduce for any
$X,Z,W\in\Gamma(\mathcal{D}^\perp)$:
\begin{eqnarray}\label{21}
       g(\overline{\nabla}_XZ,J_\alpha
W)&=&g(B(X,Z),J_\alpha W).
       \end{eqnarray}

But $M$ is a totally real ruled paraquaternionic CR-submanifold iff
$\overline{\nabla}_X Z\in\Gamma(\mathcal{D}^\perp)$, $\forall
X,Z\in\Gamma(\mathcal{D}^\perp)$ and by using (\ref{19}), (\ref{20})
and (\ref{21}) we deduce the equivalence.
\\
ii. $\Leftrightarrow$ iv. This equivalence follows from(\ref{6}).
\end{proof}

\begin{cor}\label{4.6}
Let $N$ be a paraquaternionic CR-submanifold of a paraquaternionic
K\"{a}hler manifold $(M,\sigma,g)$. If $N$ is totally geodesic, then
$N$ is a totally real ruled paraquaternionic CR-submanifold.
\end{cor}

\section{Semi-Riemannian submersions from paraquaternionic CR-submanifolds}

Semi-Riemannian submersions were introduced by O'Neill \cite{ON2}.
Let $(M,g)$ and $(M',g')$ be two connected semi-Riemannian manifolds
of index $s$ ($0\leq s\leq dim M$) and $s'$ ($0\leq s'\leq dim M'$)
respectively, with $s'\leq s$. Roughly speaking, a semi-Riemannian
submersion is a smooth map $\pi:M\rightarrow M'$ which is onto and
satisfies the following conditions:

(i) $\pi_{*|p}: T_pM \rightarrow T_{\pi(p)}M'$ is onto for all $p\in
M$;

(ii) The fibres $\pi^{-1}(p'),\ p'\in M'$, are semi-Riemannian
submanifolds of $M$;

(iii) $\pi_*$ preserves scalar products of vectors normal to fibres.

The vectors tangent to fibres are called vertical and those normal
to fibres are called horizontal. We denote by $\mathcal{V}$ the
vertical distribution, by $\mathcal{H}$ the horizontal distribution
and by $v$ and $h$ the vertical and horizontal projection. An
horizontal vector field $X$ on $M$ is said to be basic if $X$ is
$\pi$-related to a vector field $X'$ on $M'$. It is clear that every
vector field $X'$ on $M'$ has a unique horizontal lift $X$ to $M$
and $X$ is basic.

\begin{rem}\label{5.2}
If $\pi:(M,g)\rightarrow (M',g')$ is a semi-Riemannian submersion
and $X,Y$ are basic vector fields on $M$, $\pi$-related to $X'$ and
$Y'$ on $M'$, then we have the next properties (see \cite{ON}):

(i) $g(X,Y)=g'(X',Y')\circ\pi$;

(ii) $h[X,Y]$ is a basic vector field and
$\pi_*h[X,Y]=[X',Y']\circ\pi$;

(iii) $h(\nabla_XY)$ is a basic vector field $\pi$-related to
$\nabla'_{X'}Y'$, where $\nabla$ and $\nabla'$ are the Levi-Civita
connections on $M$ and $M'$;

(iv) $[E,U]\in\Gamma(\mathcal{V}), \forall U\in\Gamma(\mathcal{V})$
and $\forall E\in\Gamma(TM)$.
\end{rem}

\begin{rem}
A semi-Riemannian submersion $\pi:M\rightarrow M'$ determines, as
well as in the Riemannian case (see \cite{FIP}), two (1,2) tensor
field $T$ and $A$ on $M$, by the formulas:
\begin{equation}\label{34}
       T(E,F)=h\nabla_{vE}vF+v\nabla_{vE}hF
       \end{equation}
and respectively:
\begin{equation}\label{35}
       A(E,F)=v\nabla_{hE}hF+h\nabla_{hE}vF
       \end{equation}
for any $E,F\in \Gamma(TM)$. We remark that for
$U,V\in\Gamma(\mathcal{V})$, $T(U,V)$ coincides with the second
fundamental form of the immersion of the fibre submanifolds and for
$X,Y\in\Gamma(\mathcal{H})$, $A(X,Y)=\frac{1}{2}v[X,Y]$
characterizes the complete integrability of the horizontal
distribution $\mathcal{H}$.

It is easy to see that $T$ and $A$ satisfy:
\begin{equation}\label{35.1}
T(U,V)=T(V,U),
\end{equation}
\begin{equation}\label{35.2}
A(X,Y)=-A(Y,X)
\end{equation}
\begin{equation}\label{35.3}
\nabla_XY=h\nabla_XY+A(X,Y),
\end{equation}
\begin{equation}\label{35.4}
\nabla_UV=v\nabla_UV+T(U,V),
\end{equation}
for any $X,Y\in\Gamma(\mathcal{H})$ and $U,V\in\Gamma(\mathcal{V})$.
\end{rem}

\begin{rem}
In \cite{KOB}, S. Kobayashi observed the next similarity between a
Riemannian submersion and a CR-submanifold of a K\"{a}hler manifold:
both involve two distributions (the vertical and horizontal
distribution), one of them being integrable. Then he introduced the
concept of CR-submersion, as a Riemannian submersion from a
CR-submanifold to an almost hermitian manifold. Next, we'll consider
CR-submersions from paraquaternionic CR-submanifolds of a
paraquaternionic K\"{a}hler manifold.
\end{rem}

\begin{defn}
Let $N$ be a paraquaternionic CR-submanifold of an almost
paraquaternionic hermitian manifold $(M,\sigma,g)$ and
$(M',\sigma',g')$ be an almost \hyphenation{para-quaternionic}
hermitian manifold. A semi-Riemannian submersion $\pi:N\rightarrow
M'$ is said to be a \emph{paraquaternionic CR-submersion} if the
following conditions are satisfied:

i. $\mathcal{V}=\mathcal{D}^\perp$;

ii. For each $p\in N$, $\pi_*:\mathcal{D}_p\rightarrow T_{\pi(p)}M'$
is an isometry with respect to each complex and product structure of
$\mathcal{D}_p$ and $T_{\pi(p)}M'$, where $T_{\pi(p)}M'$ denotes the
tangent space to $M'$ at $\pi(p)$.
\end{defn}

\begin{prop}
Let $N$ be a paraquaternionic CR-submanifold of a paraquaternionic
K\"{a}hler manifold $(M,\sigma,g)$ and $(M',\sigma',g')$ be an
almost paraquaternionic hermitian manifold. If $\pi:N\rightarrow M'$
is a paraquaternionic CR-submersion, then
\begin{equation}\label{37A}
h\nabla_XJ_\alpha Y=J_\alpha
h\nabla_XY-\epsilon_{\alpha}[\omega_{\alpha+2}(X)J_{\alpha+1}Y-
    \omega_{\alpha+1}(X)J_{\alpha+2}Y]
\end{equation}
\begin{equation}\label{37}
A(X,J_\alpha Y)=J_\alpha \overline{v}B(X,Y)
\end{equation}
\begin{equation}\label{37C}
\overline{v}B(X,J_\alpha Y)=J_\alpha A(X,Y)
\end{equation}
\begin{equation}\label{37D}
\overline{h}B(X,J_\alpha Y)=J_\alpha\overline{h}B(X,Y)
\end{equation}
for any $X,Y\in\Gamma(\mathcal{D})$, where $\nabla$ is the
Levi-Civita connection on $N$ and $\overline{h}$, $\overline{v}$
denote the canonical projections on $\mathcal{\nu}$ and
$\mathcal{\nu}^\perp$, respectively.
\end{prop}
\begin{proof}
From Gauss equations and (\ref{35.3}), we have
\begin{equation}\label{36}
\overline{\nabla}_XY=h\nabla_XY+A(X,Y)+\overline{h}B(X,Y)+\overline{v}B(X,Y),
\end{equation}
for any $X,Y\in\Gamma(\mathcal{D})$, where $\overline{\nabla}$ is
the Levi-Civita connection on $M$.

From (\ref{36}) we obtain:
\begin{eqnarray}\label{37B}
(\overline{\nabla}_XJ_\alpha)Y&=&h\nabla_XJ_\alpha Y+A(X,J_\alpha Y)+\overline{h}B(X,J_\alpha Y)+\overline{v}B(X,J_\alpha Y)\nonumber\\
&&-J_\alpha h\nabla_XY-J_\alpha
A(X,Y)-J_\alpha\overline{h}B(X,Y)-J_\alpha\overline{v}B(X,Y).
\end{eqnarray}

On another hand, from (\ref{3}) we have:
\begin{equation}\label{38}
(\overline{\nabla}_XJ_{\alpha})Y=-\epsilon_{\alpha}[\omega_{\alpha+2}(X)J_{\alpha+1}Y-
\omega_{\alpha+1}(X)J_{\alpha+2}Y]\in\Gamma(\mathcal{D}),
\end{equation}
for any $X,Y\in\Gamma(\mathcal{D})$.

Comparing now (\ref{37B}) and (\ref{38}) and identifying the
components from $\Gamma(\mathcal{D})$, $\Gamma(\mathcal{D}^\perp)$,
$\Gamma(\mathcal{\nu})$ and $\Gamma(\mathcal{\nu}^\perp)$, we obtain
the wanted identities.
\end{proof}

\begin{prop}
Let $N$ be a paraquaternionic CR-submanifold of a paraquaternionic
K\"{a}hler manifold $(M,\sigma,g)$ and $(M',\sigma',g')$ be an
almost paraquaternionic hermitian manifold. If $\pi:N\rightarrow M'$
is a paraquaternionic CR-submersion, then
\begin{equation}\label{37A.2}
\overline{v}\nabla^\perp_UJ_\alpha V=J_\alpha
v\nabla_UV-\epsilon_{\alpha}[\omega_{\alpha+2}(U)J_{\alpha+1}V-
    \omega_{\alpha+1}(U)J_{\alpha+2}V]
\end{equation}
\begin{equation}\label{37.2}
J_\alpha T(U,V)=-hA_{J_\alpha V}U
\end{equation}
\begin{equation}\label{37C.2}
J_\alpha\overline{v}B(U,V)=-vA_{J_\alpha V}U
\end{equation}
\begin{equation}\label{37D.2}
\overline{h}\nabla^\perp_U J_\alpha V=J_\alpha\overline{h}B(U,V)
\end{equation}
for any $U,V\in\Gamma(\mathcal{D}^\perp)$, where $\nabla$ is the
Levi-Civita connection on $N$ and $\nabla^\perp$ is the normal
connection.
\end{prop}
\begin{proof}
From Gauss equations and (\ref{35.4}), we have:
\begin{equation}\label{36.2}
\overline{\nabla}_UV=T(U,V)+v\nabla_UV+\overline{h}B(U,V)+\overline{v}B(U,V),
\end{equation}
for any $U,V\in\Gamma(\mathcal{D}^\perp)$, where $\overline{\nabla}$
is the Levi-Civita connection on $M$.

On the other hand, from Weingarten formula, we have:
\begin{equation}\label{36.3}
\overline{\nabla}_UJ_\alpha V=-hA_{J_\alpha V}U-vA_{J_\alpha
V}U+\overline{h}\nabla^\perp_UJ_\alpha
V+\overline{v}\nabla^\perp_UJ_\alpha V,
\end{equation}
for any $U,V\in\Gamma(\mathcal{D}^\perp)$.

From (\ref{36.2}) and (\ref{36.3}) we deduce:
\begin{eqnarray}\label{37B.2}
(\overline{\nabla}_UJ_\alpha)V&=&-hA_{J_\alpha V}U-vA_{J_\alpha
V}U+\overline{h}\nabla^\perp_UJ_\alpha
V+\overline{v}\nabla^\perp_UJ_\alpha V\nonumber\\
&&-J_\alpha
T(U,V)-J_\alpha
v\nabla_UV-J_\alpha\overline{h}B(U,V)-J_\alpha\overline{v}B(U,V),
\end{eqnarray}

On another hand, from (\ref{3}) we have:
\begin{equation}\label{38.2}
(\overline{\nabla}_UJ_{\alpha})V=-\epsilon_{\alpha}[\omega_{\alpha+2}(U)J_{\alpha+1}V-
\omega_{\alpha+1}(U)J_{\alpha+2}V]\in\Gamma(\mu^\perp),
\end{equation}
for any $U,V\in\Gamma(\mathcal{D}^\perp)$.

Comparing now (\ref{37B.2}) and (\ref{38.2}) and identifying the
components from $\Gamma(\mathcal{D})$, $\Gamma(\mathcal{D}^\perp)$,
$\Gamma(\mathcal{\nu})$ and $\Gamma(\mathcal{\nu}^\perp)$, we obtain
the wanted identities.
\end{proof}

\begin{thm}
Let $N$ be a paraquaternionic CR-submanifold of a paraquaternionic
K\"{a}hler manifold  $(M,\sigma,g)$ and $(M',\sigma',g')$ be an
almost paraquaternionic hermitian manifold. If $\pi:N\rightarrow M'$
is a paraquaternionic CR-submersion, then $(M',\sigma',g')$ is a
paraquaternionic K\"{a}hler manifold.
\end{thm}
\begin{proof}
We have from (\ref{37A}):
\begin{equation}\label{39}
h\nabla_XJ_\alpha Y-J_\alpha
h\nabla_XY=-\epsilon_{\alpha}[\omega_{\alpha+2}(X)J_{\alpha+1}Y-
    \omega_{\alpha+1}(X)J_{\alpha+2}Y]
\end{equation}
for any local basic vector fields $X,Y$ on $N$.

We define the local 1-forms $\omega'_\alpha$ on $M'$ by:
\begin{equation}\label{40}
\omega'_\alpha(X')\circ\pi=\omega_\alpha(X)
\end{equation}
for any local vector field $X'$ on $M'$, $\pi$-related with an
horizontal vector field $X$ on $N$.

On the other hand, from the definition of a paraquaternionic
CR-submersion, we deduce that for any local bases
$\lbrace{J'_1,J'_2,J'_3}\rbrace$ of $\sigma'$ we have a
corresponding local basis $\lbrace{J_1,J_2,J_3}\rbrace$ of $\sigma$
such that:
\begin{equation}\label{41}
\pi_*\circ J_\alpha=J'_\alpha\circ\pi_*,\ \alpha\in\{1,2,3\}.
\end{equation}

Using (\ref{39}), (\ref{40}), (\ref{41}) and Remark \ref{5.2} we
obtain:
\begin{equation}
(\nabla'_XJ'_{\alpha})Y=-\epsilon_{\alpha}[\omega'_{\alpha+2}(X')J'_{\alpha+1}Y'-
\omega'_{\alpha+1}(X')J'_{\alpha+2}Y']
\end{equation}
for any local vector fields $X',Y'$ on $M'$, $\pi$-related with two
local basic vector fields $X,Y$ on $N$, where $\nabla'$ is the
Levi-Civita connection on $M'$. Hence $(M',\sigma',g')$  is a
paraquaternionic K\"{a}hler manifold.
\end{proof}

\begin{thm}
Let $N$ be a mixed foliated paraquaternionic CR-submanifold of a
paraquaternionic K\"{a}hler manifold $(M,\sigma,g)$ and
$(M',\sigma',g')$ be an almost paraquaternionic hermitian manifold.
If $\pi:N\rightarrow M'$ is a paraquaternionic CR-submersion, then
 $N$ is locally a semi-Riemannian product of a paraquaternionic
submanifold and a totally real submanifold of $M$. In particular, if
$N$ is complete and simply connected then it is a global
semi-Riemannian product.
\end{thm}
\begin{proof}
Because $N$ is a mixed foliated paraquaternionic CR-submanifold, in
particular $N$ is a mixed geodesic submanifold, and from (\ref{18})
we obtain that the distribution $\mathcal{D}^\perp$ is parallel. On
another hand, $\mathcal{D}$ being integrable, we have
$A(X,Y)=v\nabla_XY=0,$ for any $X,Y\in\Gamma(\mathcal{D})$ and
therefore the distribution $\mathcal{D}^\perp$ is also parallel.
Hence $N$ is locally a semi-Riemannian product
$(N_1,g_1)\times(N_2,g_2)$, where $N_1$ and $N_2$ are leaves of
$\mathcal{D}$ and $\mathcal{D}^\perp$.

Finally, if $N$ is complete and simply connected, applying the
decomposition theorem for semi-Riemannian manifolds (see \cite{WU})
we obtain the last part of the theorem.
\end{proof}

\begin{ex}
We remarked in Section 3 that \[ N=\{[u_0,u_1,u_2]\in
f_{p,q}^{-1}(0)|q^2|u_0|^2+p^2|u_1|^2+p^2|u_2|^2\neq 0\}
\]
is a a proper paraquaternionic $CR$-submanifold of $P^2\mathbb{B}$.
Moreover, the Lie group $G=\{e^{jt}|t\in\mathbb{R}\}$ acts freely
and isometrically on $N$. Using now the paraquaternionic K\"{a}hler
reduction (see Theorem 5.2 from \cite{VKM}) we obtain that the
manifold $M'=N/G$ equipped with the submersed metric (i.e the one
$g'$ which makes the projection $\pi:(N,g)\rightarrow (M',g')$ a
semi-Riemannian submersion) is again a paraquaternionic K\"{a}hler
manifold with respect to the structure $\sigma'$ induced on $M'$
from the structure $\sigma$ by the projection $\pi$. Moreover,
$\pi:N\rightarrow N/G$ is a paraquaternionic CR-submersion.
\end{ex}

\section{Some curvature properties of paraquaternionic CR-submersions}

Let $(M,g)$ be a semi-Riemannian manifold. The sectional curvature
$K$ of a 2-plane in $T_pM,$ $p\in M$, spanned by
$\lbrace{X,Y}\rbrace$, is defined by:
     \begin{equation}\label{1.def}
          K(X\wedge Y)=\frac{R(X,Y,X,Y)}{g(X,X)g(Y,Y)-g(X,Y)^2}.
    \end{equation}
It is clear that the above definition makes sense only for
non-degenerate planes, i.e. those satisfying $Q(X\wedge
Y)=g(X,X)g(Y,Y)-g(X,Y)^2\neq 0.$

If $\pi:(N,g)\rightarrow (M',g')$ is a semi-Riemannian submersion,
then the tensor fields $A$ and $T$ defined in the above Section play
a fundamental role in expressing the curvatures of the $(N,g)$,
$(M',g')$ and of any fibre $(\pi^{-1}(p'),\hat{g}_{p'})$, $p'\in
M'$, because we have the following formulas stated in \cite{ON}:
\begin{eqnarray}\label{42}
R(U,V,W,W')&=&\hat{R}(U,V,W,W')-g(T(U,W),T(V,W'))\nonumber\\&&+g(T(V,W),T(U,W')),
\end{eqnarray}
\begin{eqnarray}\label{42.1}
R(X,Y,Z,Z')&=&R^*(X,Y,Z,Z')-2g(A(X,Y),A(Z,Z'))\nonumber\\&&+g(A(Y,Z),A(X,Z'))-g(A(X,Z),A(Y,Z')),
\end{eqnarray}
for any $U,V,W,W'\in\Gamma(\mathcal{V})$ and
$X,Y,Z,Z'\in\Gamma(\mathcal{H})$, where $R$ is the Riemannian
curvature of $(N,g)$, $\hat{R}$ is the Riemannian curvature of fibre
and $R^*(X,Y,Z,Z')=g(R^*(Z,Z')Y,X)$, $R^*(\cdot,\cdot)\cdot$ being
the (1,3)-tensor field on $\Gamma(\mathcal{H})$ with values in
$\Gamma(\mathcal{H})$ which associates to any
$X,Y,Z\in\Gamma(\mathcal{H})$ and $p\in N$, the horizontal lift
$(R^*(X,Y)Z)_p$ of
$R'_{\pi(p)}(\pi_{*p}(X_p),\pi_{*p}(Y_p))\pi_{*p}(Z_p)$, $R'$
denoting the Riemannian curvature of $(M',g')$.

\begin{thm}
Let $N$ be a paraquaternionic CR-submanifold of a paraquaternionic
K\"{a}hler manifold $(M,\sigma,g)$ and $(M',\sigma',g')$ be an
almost paraquaternionic hermitian manifold. If $\pi:N\rightarrow M'$
is a paraquaternionic CR-submersion, then the sectional curvatures
of $M$ and the fibres are related by:
\begin{eqnarray}\label{1.AA}
\overline{K}(U\wedge V)=\hat{K}(U\wedge
V)&-&\epsilon_\alpha\theta_U\theta_V[g(A_{J_\alpha
U}U,A_{J_\alpha V}V)-g(A_{J_\alpha V}U,A_{J_\alpha V}U)]\nonumber\\
&-&\epsilon_\alpha\theta_U\theta_V[g(\overline{h}\nabla^\perp_UJ_\alpha
U,\overline{h}\nabla^\perp_VJ_\alpha
V)-g(\overline{h}\nabla^\perp_UJ_\alpha
V,\overline{h}\nabla^\perp_UJ_\alpha V)]
\end{eqnarray}
for any unit space-like or time-like orthogonally vector fields
$U,V\in\Gamma(\mathcal{D}^\perp)$ and $\alpha\in\{1,2,3\}$, where
$\theta_U=g(U,U)\in\{-1,1\}$ and $\theta_V=g(V,V)\in\{-1,1\}$.
\end{thm}
\begin{proof}
From (\ref{42}) and Gauss equation we have:
\begin{eqnarray}\label{43}
\overline{R}(U,V,W,W')&=&\hat{R}(U,V,W,W')-g(T(U,W),T(V,W'))+g(T(V,W),T(U,W'))\nonumber\\
&&-g(B(U,W),B(V,W'))+g(B(V,W),B(U,W')),
\end{eqnarray}
for any $U,V,W,W'\in\Gamma(\mathcal{D}^\perp)$.

From (\ref{1.def}) we see that the sectional curvatures of $M$ and
the fibres for a 2-plane spanned by two unit space-like or time-like
orthogonally vector fields $U,V$ are:
\begin{equation}\label{44}
\overline K(U\wedge V)=\theta_U\theta_V\overline{R}(U,V,U,V)
\end{equation}
and, respectively,
\begin{equation}\label{45}
\hat{K}(U\wedge V)=\theta_U\theta_V\hat{R}(U,V,U,V),
\end{equation}
where $\theta_U=g(U,U)\in\{-1,1\}$ and $\theta_V=g(V,V)\in\{-1,1\}$.

From (\ref{43}), (\ref{44}) and (\ref{45}), using (\ref{35.1}), we
obtain:
\begin{eqnarray}\label{46}
\overline{K}(U\wedge V)&=&\hat{K}(U\wedge V)-\theta_U\theta_V[g(T(U,U),T(V,V))-g(T(U,V),T(U,V))\nonumber\\
&&+g(B(U,U),B(V,V))-g(B(U,V),B(U,V))]
\end{eqnarray}
for any unit space-like or time-like orthogonally vector fields
$U,V\in\Gamma(\mathcal{D}^\perp)$.

Finally, using (\ref{37.2}), (\ref{37C.2}) and (\ref{37D.2}) in
(\ref{46}), we obtain (\ref{1.AA}).
\end{proof}

\begin{thm}
Let $N$ be a paraquaternionic CR-submanifold of a paraquaternionic
K\"{a}hler manifold $(M,\sigma,g)$ and $(M',\sigma',g')$ be an
almost paraquaternionic hermitian manifold. If $\pi:N\rightarrow M'$
is a paraquaternionic CR-submersion, then for any unit space-like or
time-like horizontal vector field $X$ one has:
\begin{equation}\label{2.AA}
\overline{H}_{\alpha}(X)=H'_{\alpha}(\pi_*X)-4g(\overline{v}B(X,X),\overline{v}B(X,X))+
2g(\overline{h}B(X,X),\overline{h}B(X,X)),
\end{equation}
for $\alpha\in\{1,2,3\}$, where $\overline{H}_{\alpha}$ and
$H'_{\alpha}$ are the holomorphic sectional curvatures of $M$ and
$M'$, defined by $\overline{H}_{\alpha}(X)=\overline{K}(X\wedge
J_\alpha X)$ and $H'_{\alpha}(X)=K'(X\wedge J_\alpha X)$,
respectively.
\end{thm}
\begin{proof}
From (\ref{42.1}) and Gauss equation we have:
\begin{eqnarray}\label{100}
\overline{R}(X,Y,Z,Z')&=&R^*(X,Y,Z,Z')-2g(A(X,Y),A(Z,Z'))\nonumber\\
&&+g(A(Y,Z),A(X,Z'))-g(A(X,Z),A(Y,Z'))\nonumber\\
&&+g(B(Y,Z),B(X,Z'))-g(B(Y,Z'),B(X,Z)),
\end{eqnarray}
for any $X,Y,Z,Z'\in\Gamma(\mathcal{D})$.

From (\ref{2}) and (\ref{1.def}) we see that the holomorphic
sectional curvature of $M$ is:
\begin{equation}\label{101}
\overline{H}_{\alpha}(X)=\epsilon_\alpha\overline{R}(X,J_\alpha
X,X,J_\alpha X)
\end{equation}
for any unit space-like or time-like vector field
$X\in\Gamma(\mathcal{D})$.

On the other hand we can easily see that:
\begin{equation}\label{102}
H'_{\alpha}(\pi_*X)=\epsilon_\alpha R^*(X,J_\alpha X,X,J_\alpha X)
\end{equation}
for any unit space-like or time-like vector field
$X\in\Gamma(\mathcal{D})$.

From (\ref{100}), (\ref{101}) and (\ref{102}), using (\ref{35.2}),
we obtain:
\begin{eqnarray}\label{103}
\overline{H}_\alpha(X)&=&H'_\alpha(\pi_*X)-3\epsilon_\alpha g(A(X,J_\alpha X),A(X,J_\alpha X))\nonumber\\
&&+\epsilon_\alpha[g(B(X,J_\alpha X),B(X,J_\alpha
X))-g(B(X,X),B(J_\alpha X,J_\alpha X))].
\end{eqnarray}

From (\ref{35.2}) and (\ref{37C}) one has:
\begin{equation}\nonumber
0=A(X,X)=-\epsilon_\alpha J_\alpha\overline{v}B(X,J_\alpha X)
\end{equation}
and so we derive:
\begin{equation}\label{104}
\overline{v}B(X,J_\alpha X)=0.
\end{equation}

Similarly, from (\ref{35.2}), (\ref{37}) and (\ref{37C}) we derive:
\begin{equation}\label{105}
\overline{v}B(J_\alpha X,J_\alpha X)=\epsilon_\alpha
\overline{v}B(X,X).
\end{equation}
and from (\ref{37D}), using the symmetry of $B$, we obtain:
\begin{equation}\label{106}
\overline{h}B(J_\alpha X,J_\alpha X)=-\epsilon_\alpha
\overline{h}B(X,X).
\end{equation}

Finally, using (\ref{104}), (\ref{105}) and (\ref{106}) in
(\ref{103}) we obtain (\ref{2.AA}).
\end{proof}

\begin{cor}
Let $N$ be a totally geodesic paraquaternionic CR-submanifold of a
paraquaternionic K\"{a}hler manifold $(M,\sigma,g)$ and
$(M',\sigma',g')$ be an almost paraquaternionic hermitian manifold.
If $\pi:N\rightarrow M'$ is a paraquaternionic CR-submersion one
has:
\begin{equation}\label{2.AAA}
\overline{H}_{\alpha}(X)=H'_{\alpha}(\pi_*X),
\end{equation}
for any unit space-like or time-like horizontal vector field $X$.
\end{cor}

\section*{Acknowledgments}
We would like to express our sincere thanks to Professor Liviu Ornea
for valuable discussions concerning this subject. S. Ianu\c{s} and
G.E. V\^{\i}lcu are supported by a CNCSIS PN2-IDEI grant, no.
525/2009. S. Marchiafava activity was done under the programs of
G.N.S.A.G.A. of C.N.R.-Indam and "Riemannian metrics and
differentiable structures" of MIUR (Italy).

\newpage

\begin{center}
Stere Ianu\c s \\
{\em    University of Bucharest,
        Department of Mathematics,\\
        Str. Academiei, Nr. 14, Bucharest 70109, Romania}\\
        e-mail: ianus@gta.math.unibuc.ro\\
\end{center}

\begin{center}
Stefano Marchiafava \\
{\em    Universit\`{a} di Roma "La Sapienza",
        Dipartimento di Matematica,\\
        P.le Aldo Moro N. 2, 00185 Rome, Italy}\\
        e-mail: marchiafava@mat.uniroma1.it\\
\end{center}

\begin{center}
Gabriel Eduard V\^\i lcu$^{1,2}$ \\
{\em
      $^1$University of Bucharest, Faculty of Mathematics and Computer Science,\\
      Research Center in Geometry, Topology and Algebra,\\
      Str. Academiei, Nr. 14, Sector 1, Bucharest 70109, Romania\\
      $^2$Petroleum-Gas University of Ploie\c sti,\\
      Department of Mathematics and Computer Science,\\
      Bulevardul Bucure\c sti, Nr. 39, Ploie\c sti 100680, Romania}\\
      e-mail: gvilcu@mail.upg-ploiesti.ro
\end{center}

\end{document}